\documentclass{paper}

\usepackage{amsmath}
\usepackage{verbatim}
\usepackage{amssymb}
\usepackage{latexsym}
\usepackage[english]{babel}
\usepackage{amsthm,dsfont,mathrsfs}
\usepackage[T1]{fontenc}

\usepackage[usenames,dvipsnames]{pstricks}
\usepackage{pstricks-add,pst-plot,pst-node}
\usepackage{epsfig}
\usepackage{graphicx}

\usepackage{enumerate}
\newtheorem{Theorem}{Theorem}
\newtheorem{Corollary}[Theorem]{Corollary}

\newtheorem{Proposition}[Theorem]{Proposition}

\newtheorem{Assumption}[Theorem]{Assumption}

\newenvironment{Proof}{\vspace{1ex}{\sc Proof}.}{\vspace{1ex}}

\def\Re{\mathds{R}}

\def\1e{\mathds{1}}

\def\vg0{\mathbf{0}}

\def\A{\mathscr{A}}

\def\D{\mathscr{D}}

\def\I{\mathscr{I}}

\def\L{\mathscr{L}}

\def\P{\mathscr{P}}
\def\Q{\mathscr{Q}}

\def\e{\varepsilon}

\def\eop{\rule{0.7ex}{0.7ex}}
\def\eqd{: =}

\def\segoo#1#2{(#1, #2)}

\def\scal#1#2{\langle #1,#2 \rangle}

\def\set#1#2{\left\{\mskip 1mu #1 \mskip 1mu
    | \mskip 1mu #2 \mskip 1mu \right\}
    }

\def\norm#1{\left\| \mskip 1mu #1 \mskip 1mu \right\|}

\def\co#1{\mathop{\mathrm{co}}#1}

\def\argmin{\mathop{\mathrm{argmin}}}

\begin{document}

\parindent0pt
\parskip8pt
%-------------------
\title{{\bf
Lagrange duality for the Morozov principle}}
%A note on the Morozov principle via Lagrange duality
%Regularization of linear ill-posed problems: from noise assessment to parameter estimation}}
\author{Xavier Bonnefond\thanks{MAPMO, Université d'Orléans, France} and Pierre Mar\'echal\thanks{ISAE, Université de Toulouse, France}}
\date{\today}
%-------------------
\maketitle

\abstract
Considering a general linear ill-posed equation, we explore the duality arising from the requirement that the discrepancy should take a given value based on the estimation of the noise level, as is notably the case when using the Morozov principle. We show that, under reasonable assumptions, the dual function is smooth, and that
its maximization points out the appropriate value of Tikhonov's regularization parameter.

%%%%%%%%%%%%%%%%%%%%%%
\section{Introduction}

Let us consider the ill-posed linear inverse problem
\begin{equation}
\label{Af-eg-g}
\A f=g,
\end{equation}
in which $\A$ maps the Hilbert space $F$ into the Hilbert space $G$
linearly, $g\in G$ is the data and $f\in F$ is the unknown. As usual,
we assume that
$$
g=g_0+\delta g
$$
where $g_0\eqd \A f_0$ for some $f_0\in F$ and
$\delta g$ is the unknown noise. In many applications,
the ill-posedness of~\eqref{Af-eg-g} is the consequence of the
compactness of~$\A$. Numerous strategies have been developped in
the last decades to regularize such problems. The variational
approach consists in defining the reconstructed object as the
minimizer of some functional, which usually is the sum of a
\emph{fit term} and of a \emph{regularization term}.
More precisely, in this approach, a family
of such functionals is considered, which depends on one parameter~$\alpha$
(or more). A natural requirement is that, for positive values of
$\alpha$, the correponding variational problem is well-posed, while
letting $\alpha\downarrow 0$ yields, at the limit, a least
square solution of Problem~\eqref{Af-eg-g}.

In practice, the choice of $\alpha$ is a crucial step. As matter
of fact, large values of~$\alpha$ correspond to coarse approximations
of the original model, while small values cause high sensitivity
of the solution to perturbations on the data side (which we may
call \emph{hypersensitivity}).

Strategies for the determination of~$\alpha$ may be cast into \emph{a priori}
and \emph{a posteriori} approaches. An example of the first class consists
in selecting~$\alpha$ so as to satisfy some stability requirement, regardless
of the particular data~$g$ to be considered. Such a choice may
yield a value of~$\alpha$ which \emph{overregularizes} the problem.
This is why one usually focuses on \emph{a posteriori}
parameter selection strategies, that is, strategies that are
\emph{data dependent}.

It is worth mentioning that, in~\cite{iveta}, it was shown that the Golub-Kahan
bidiagonalization algorithm enables the estimation of the noise level.
This was achieved by observing the corruption by noise of the iterates produced by the
algorithm.
Of course, an estimation of the noise level in the data may also be obtained
from the modeling of the data acquisition process.
At all events,
it then seems reasonable to find a value of~$\alpha$ such that
the corresponding solution~$f_\alpha$ of Problem~$(\P_\alpha)$
satisfies
\begin{equation}
\label{NoiseConstraint}
\norm{\A f_\alpha-g}^2\simeq\tau^2,
\end{equation}
in which~$\tau$ is an estimation of $\norm{\delta g}$.
The Morozov Principle states that $\alpha$ should be chosen
in such a way that Eq.~\eqref{NoiseConstraint} is satisfied
exactly with $\tau$ replaced by $c\tau$, where $c$ is a constant
strictly greater than~1, and in fact close to~1. This principle
is often used as a stopping criterion for iterative regularization
schemes (see, e.g, \cite{frick}).

In this note, we show how the estimation of~$\tau$ (whatever may be the estimation
method) can be used to directly determine~$\alpha$.
This will be done by \emph{dualizing} the constraint on the noise level.

In Section~\ref{notation}, we fix the context and make comments
on the interplay between the various manners to relax the
initial constraint equation~\eqref{Af-eg-g}.
In Section~\ref{duality-for-morozov}, we explore the duality arising
from a constraint of the form~\eqref{NoiseConstraint}.

A smooth reading of the next sections may require some familiarity with
variational and convex analysis. Our reference books in convex analysis
are~\cite{rockafellar,hull-1,hull-2,Zalinescu}.

%%%%%%%%%%%%%%%%%%%%%%%%%%%%%%%%%%%%%%%%%%
\section{Notation and preliminary remarks}
\label{notation}

The spaces~$F$ and~$G$ are endowed with the
norms $\norm{\cdot}_F$ and $\norm{\cdot}_G$ 
associated with the inner products $\scal{\cdot}{\cdot}_F$ and
$\scal{\cdot}{\cdot}_G$, respectively. We shall frequently omit the
subscripts $F$ and $G$, since most of the time the context leaves no
ambiguity.

Most variational regularization techniques consist in defining 
the reconstructed object as the solution of
$$
(\P_\alpha)\quad\left|
\begin{array}{rl}
\hbox{Minimize}&\norm{\A f-g}^2+\alpha J(f)\\
\hbox{s.t.}&f\in F,
\end{array}
\right.
$$
in which~$J$ is the so called \emph{regularizer}.
It is customary to make the following assumption:
\begin{Assumption}\sf
\label{ConvexSciCoercive}
The function
$J$ is proper convex, lower semi-continuous and coercive.
\end{Assumption}
Recall that a function $J(f)$ is said to be coercive if $J(f)\to\infty$ as $\norm{f}\to\infty$.
In this paper, we also make throughout the additional (reasonable) assumption:
\begin{Assumption}\sf
\label{StrictConvexityKernel}
The function
$J$ is strictly convex along~$\ker\A$ and attains its minimum on~$\ker\A$.
\end{Assumption}
It is well known that the solution $f_\alpha$ to Problem~$(\P_\alpha)$
is also solution to the following constrained problem:
$$
(\Q)\quad
\left|
\begin{array}{rl}
\hbox{Minimize}&J(f)\\
\hbox{s.t.}&\norm{\A f-g}^2=\e
\end{array}
\right.
$$
with $\e=\e(\alpha)\eqd\norm{\A f_\alpha-g}^2$.
This is a particular case of Everett's lemma (see~\cite{hull-2}, for example).
Lagrange duality makes it possible to go the other way around: starting from
a problem such as $(\Q)$ with $\e=\tau^2$ (since we wish to prescribe the tolerance~$\tau$),
we may compute the value of $\alpha$
ensuring that $f_\alpha$ is also solution to~$(\P_\alpha)$.

Notice that Problem~$(\Q)$ is not convex, but that whenever 
\begin{equation}
\label{taille_bruit}
\norm{g}^2\geq\e
\end{equation}
(which is a natural assumption since the noise is usually smaller than the data),
Problem~$(\Q)$ is equivalent to
$$
(\Q^*)\quad
\left|
\begin{array}{rl}
\hbox{Minimize}&J(f)\\
\hbox{s.t.}&\norm{\A f-g}^2\leq\e
\end{array}
\right.
$$
Indeed the only case in which the solutions of these problems are different occurs
when the solution~$f^*$ of Problem~$(\Q^*)$ satisfies~$\norm{\A f^*-g}^2<\e$. This
means that the constraint is not active and that the optimality condition reads
$$
0\in\partial J(f^*).
$$
According to Assumption~\ref{StrictConvexityKernel}, this yields~$\A f^*=0$, so that~$\norm{g}^2<\e$,
in contradiction with~\eqref{taille_bruit}.

For convenience, we shall speak of \emph{tolerance} or \emph{penalized}
formulations in order to refer to problems such as~$(\Q^*)$ or~$(\P_\alpha)$, respectively.

We emphasize that, although solving $(\Q^*)$ with $\e=\tau^2$ directly
is possible in principle, it is not satisfactory in practice. The equivalence
between the penalized and tolerance formulations does not apply to
stability analysis. An enlightening example is provided by the case where~$J(f)$ 
is~$\norm{f}^2$ (Tikhonov regularization): in this case, the stability involves, in
the penalized formulation, the spectral
properties of $\A^*\A+\alpha\I$, while the tolerance formulation may retain initial
instability if~$\tau=\mathrm{dist}(g,\overline{\A F})$ (since then the unique solution 
is just the \emph{unstable} minimum norm least square solution).

%%%%%%%%%%%%%%%%%%%%%%%%%%%%%%%%%%%%%%%%%%%
\section{Duality for the Morozov principle}
\label{duality-for-morozov}

From now on, we fix $\e=\tau^2$ in Problem~$(\Q^*)$.
The Lagrangian of~$(\Q^*)$ is given by
$$
L(f,\lambda)\eqd J(f)+\lambda\left(\norm{\A f-g}^2-\e\right),\quad
f\in F,\;\lambda\in\Re_+,
$$
and the Lagrange problem associated to~$(\Q)$ and~$\lambda$ reads
$$
(\L_\lambda)\left|
\begin{array}{rl}
\hbox{Minimize}&L(f,\lambda)\\
\hbox{s.t.}&f\in F.
\end{array}
\right.
$$
We see right away that the above Lagrange problem is equivalent, for $\lambda>0$, to
the Tikhonov problem $(\P_\alpha)$ with $\alpha=1/\lambda$. It is convex, and the first 
order optimality condition reads:
\begin{equation}
\label{opt_L_lambda}
0\in\partial J(f)+2\lambda\left(\A^*\A f-\A^* g\right).
\end{equation}
From Assumptions~\ref{ConvexSciCoercive} and~\ref{StrictConvexityKernel}, it is readily 
seen that, for every~$\lambda>0$, Problem~$(\L_\lambda)$ has a unique solution~$f_\lambda$ 
satisfying~\eqref{opt_L_lambda}.

The \emph{dual function} is defined as the optimal value in the
Lagrange problem:
$$
D(\lambda)\eqd\inf\set{L(f,\lambda)}{f\in F}=L(f_\lambda,\lambda),\quad\lambda\in\Re,
$$
and the dual problem associated with~$(\Q)$ is:
$$
\left|
\begin{array}{rl}
\hbox{Maximize}&D(\lambda)\\
\hbox{s.t.}&\lambda\in\Re,
\end{array}
\right.
$$
which is obviously equivalent to
$$
(\D)
\left|
\begin{array}{rl}
\hbox{Maximize}&D(\lambda)\\
\hbox{s.t.}&\lambda>0.
\end{array}
\right.
$$

Before stating the main result of this section, we recall an important
result on the subdifferential of a supremum of convex functions.

\begin{Theorem}\sf
\label{subdifferential-sup-convex-functions}
Let $Y$ be a compact set in some metric space, and let $\varphi\colon Y\times\Re^n\to\Re$
be such that
\begin{enumerate}
\item[(1)]
for every $y\in Y$, $\varphi(y,\cdot)$ is convex;
\item[(2)]
for every $x\in\Re^n$, $\varphi(\cdot,x)$ is upper semicontinuous.
\end{enumerate}
If $x$ is a point where the convex function $\Phi(x)\eqd\sup_{y\in Y}\varphi(y,x)$
has a compact subdifferential ($x$ must be in the interior of the effective domain
of~$\Phi$), then
$$
\partial\Phi(x)=\co{\bigcup_{y\in Y(x)}\partial\varphi(y,\cdot)(x)},
$$
in which $Y(x)\eqd\set{y\in Y}{\varphi(y,x)=\Phi(x)}$ and $\co{S}$ denotes
the convex hull of the set~$S$.
\end{Theorem}

This classical result from subdifferential calculus has more general forms, and the interested reader may
consult \cite{Hantoute} and the references therein. We shall only need the following corollaries of the above theorem.

\begin{Corollary}\sf
\label{Cor1}
With the assumption of the theorem, suppose in addition that,
for every $y\in Y$, $\varphi(y,\cdot)$ is differentiable. Then,
$$
\partial\Phi(x)=\co\set{\nabla\varphi(y,\cdot)(x)}{y\in Y(x)}.
$$
\end{Corollary}

\begin{Corollary}\sf
\label{Cor2}
With the assumptions of the theorem and the previous corollary, suppose in addition that,
at the point~$x$, $\varphi(\cdot,x)$ attains its maximum at a unique $y=y(x)\in Y$. Then~$\Phi$
is differentiable at~$x$ and
$$
\nabla\Phi(x)=\nabla\varphi\big(y(x),\cdot\big)(x).
$$
\end{Corollary}

\begin{Proposition}\sf
The dual function $D$ is differentiable on~$\segoo{0}{\infty}$, and its
derivative is given by $D'(\lambda)=\norm{\A f_\lambda -g}^2-\e$.
\end{Proposition}

\begin{Proof}
As the infimum of a collection of affine functions, $D$ is concave.
Since, for every $\lambda>0$, Problem~$(\L_\lambda)$ has a unique solution,
we can apply Corollary~\ref{Cor2}.~\eop
\end{Proof}

\begin{Theorem}\sf
\label{estimation-lambda}
Suppose that the noise estimation~$\tau$ and the data~$g$ satisfies:
\begin{equation}
\label{ineq_tau}
\mathrm{dist}(g,\overline{\A F})<\tau<\norm{g}.
\end{equation}
Then Problem~$(\D)$ as at least one solution~$\bar\lambda>0$, and
the unique solution~$f_{\bar\lambda}$ of Problem~$(\L_{\bar\lambda})$
satisfies~$\norm{\A f_{\bar\lambda} -g}^2=\tau^2$,
and is consequently a solution of Problem~$(\Q)$ too.
Moreover, any other solution~$\bar\lambda'>0$ of Problem~$(\D)$ leads to
the same~$f_{\bar\lambda'}=f_{\bar\lambda}$.
\end{Theorem}

\begin{Proof}
Clearly, $D(0)=0$ and $D'(0)=\norm{g}^2-\e>0$, in which $D'(0)$ denotes the right derivative of~$D$ at~$0$.
In addition, we have
$
\mathrm{dist}(g,\overline{\A F})<\tau,
$
so that there exists some~$f_0\in F$ such that~$\norm{\A f_0-g}^2-\e<0$. Then,
$$
D(\lambda)\leq L(f_0,\lambda)=J(f_0)+\lambda(\norm{\A f_0-g}^2-\e)
\longrightarrow -\infty
\quad\hbox{as}\quad
\lambda\to\infty,
$$
and this is sufficient to prove that Problem~$(\D)$ has a solution~$\bar\lambda$. One has
$$
D'(\bar\lambda)=0=\norm{\A f_{\bar\lambda}-g}^2-\e,
$$
so that~$f_{\bar\lambda}$ is actually solution of Problem~$(\Q)$.
Now, let~$\bar\lambda$ and~$\bar\lambda'$ be two solutions of Problem~$(\D)$ and
let~$\bar f:=f_{\bar\lambda}\neq \bar f':=f_{\bar\lambda'}$. Note that,
since~$L(\bar f,\bar\lambda)=L(\bar f',\bar\lambda')$, we 
have~$J(\bar f)=J(\bar f')$. 
Using the optimality condition~\eqref{opt_L_lambda} at~$\bar f$ and~$\bar f'$ consecutively, 
one gets:
$$
J(\bar f')\geq J(\bar f)-\scal{2\bar\lambda\A^*(\A \bar f-g)}{\bar f'-\bar f}
$$
and
$$
J(\bar f)\geq J(\bar f')-\scal{2\bar\lambda'\A^*(\A \bar f'-g)}{\bar f-\bar f'}.
$$
This yields:
$$
0\leq\scal{2\A^*(\A \bar f-g)}{\bar f'-\bar f}
\quad\hbox{and}\quad
0\leq \scal{2\A^*(\A \bar f'-g)}{\bar f-\bar f'}.
$$
Subtracting the last two inequalities, we get
$$
\norm{\A (\bar f -\bar f')}^2\leq 0,
$$
so that~$\bar f-\bar f'\in\ker \A$. Since~$J$ satisfies
Assumption~\ref{StrictConvexityKernel}, the element
$$
\bar f'':=\frac{\bar f +\bar f'}{2}
$$
satisfies~$J(\bar f'')<J(\bar f)$ and~$\A \bar f''=\A \bar f$. Finally, we get
$
L(\bar f'',\bar\lambda)<L(\bar f,\bar\lambda),
$
which contradicts the optimality of~$\bar f$. In conclusion,~$\bar f=\bar f'$.~\eop
\end{Proof}

The theorem opens the way to many iterative algorithms for the search of
a solution of~$(\Q)$. For example:

{\bf Algorithm.}\vspace{-1ex}
\begin{enumerate}[(i)]
\item
Initialization: $\lambda_0=0$.
\item
Iteration:\\[.5ex]
$f_n=\argmin L(\cdot,\lambda_{n-1})$,\\
$\lambda_n=\lambda_{n-1}+\rho_n(D'(\lambda_{n-1}))$.
\end{enumerate}

Here $\rho_n$ is the step size, which depends on the selected method to maximize~$D$.
This algorithm is to be compared to the augmented Lagrangian method as described for example
in~\cite{frick}.

In figure~\ref{figure-aspect-D} below, we sketch the behavior of~$D$
according the accuracy of the estimation of the noise level.
The solid line corresponds to the assumption~\eqref{ineq_tau} in the above theorem,
which ensures that~$D$ has a maximum on $\segoo{0}{\infty}$.
The dotted line corresponds to the case where the data is dominated by the noise.
In the case of the dashed line, the dual function does not attain a maximum because the
estimation of the noise is too optimistic. Such a behavior would occur for example if
one takes $\varepsilon=0$, in which case the constraint in Problem~$(\Q)$ is equivalent
to the equality constraint~\eqref{Af-eg-g}.
In this last case, the above algorithm requires a stopping criterion, which usually
satisfies the Morozov discrepancy principle, that is to say, the condition
$\norm{\A f -g}\leq\tau$.

\begin{figure}[ht]
\centering
\includegraphics[width=\textwidth]{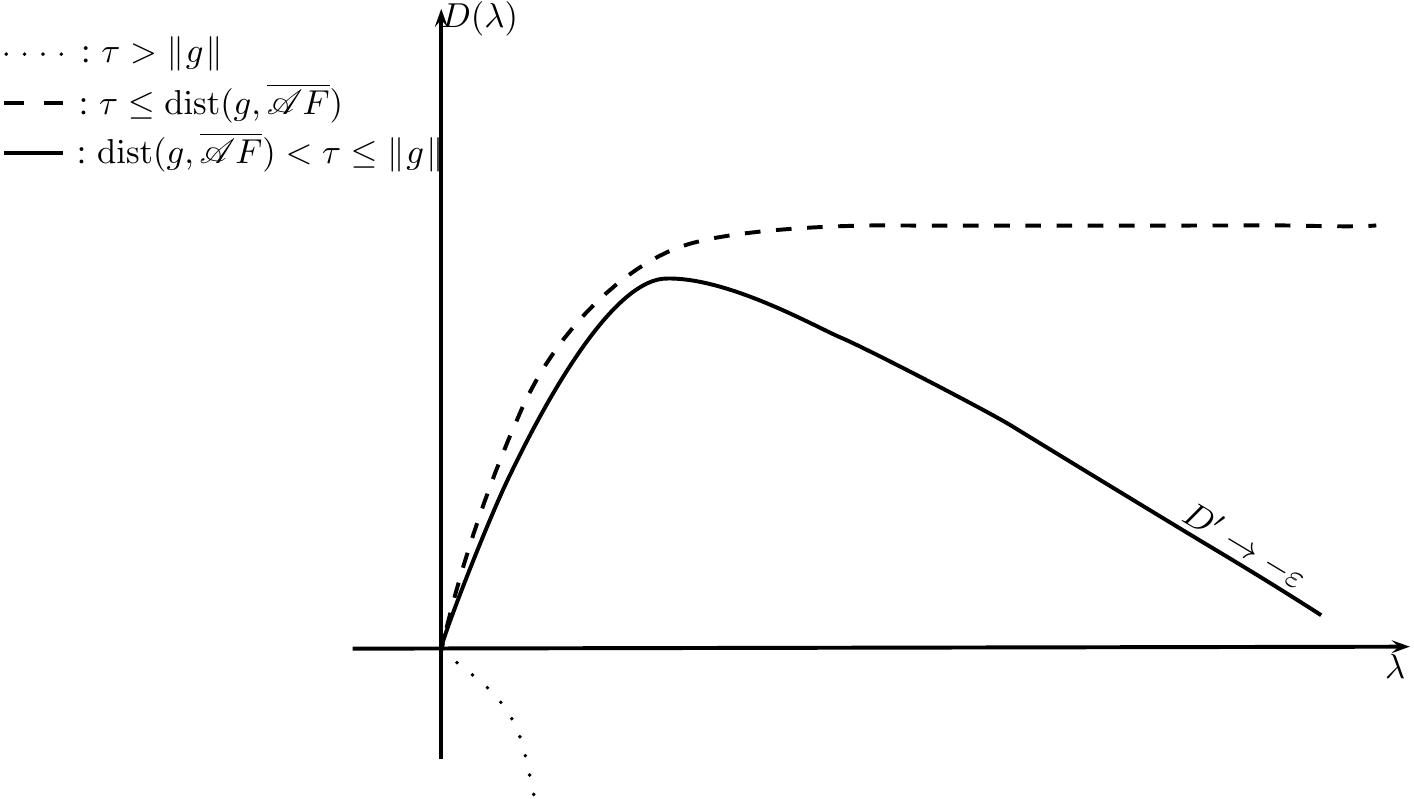}
\caption{Aspect of the dual function $D$ in several situations}
\label{figure-aspect-D}
\end{figure}

On the one hand, if the condition $\mathrm{dist}(g,\overline{\A F})<\tau<\norm{g}$
is not satisfied, the maximization of~$D$ can only fail, so it will be
clear that the desired tolerance cannot be reached. This may happen
if the noise estimation is not sufficiently accurate, or if the
noise level dominates the data.

On the other hand, if the estimation of~$\tau$ is reasonably
accurate, the condition
$$
\mathrm{dist}(g,\overline{\A F})<\tau<\norm{g}
$$
is likely to be satisfied, and the desired~$\alpha$ can be computed
by maximizing $D(\lambda)$. The evaluation of~$D$ and~$D'$ can be performed by
solving a Lagrange problem. Such
problems are well behaved for small values of~$\lambda$, and their
condition will deteriorate as~$\lambda$ will increase.

\section{Conclusion}
In this note, we have explored some aspects of the dualization
of the constraint arising from the implementation of the Morozov
principle. We have shown that, under mild assumptions, the dual
function is smooth, and that its shape is directly related to the
magnitude of the noise and the quality of its estimation. In the
favorable cases, the desired Tikhonov parameter can be obtained easily
via the maximization of the dual function.

\bibliographystyle{plain}
\bibliography{tikho-mollif}

\end{document}